\documentstyle{article}

\def \E {{\bf E\, }}

\def \l {\lambda}

\def \g {\gamma}

\def \t {\tilde}

\def \vep {\varepsilon}

\def \b {\beta}

\def \a {\alpha}

\def \d {\delta}

\def \N {{\cal N}}

\def \L {{\cal L}}

\def \fb {f^{(\b)}}

\def \Z {\hat Z_{N,m}(\beta)}

\def \I {{\bf I}}

\def \P {{\hbox{Pr}}}

\def \lf {\lfloor}

\def \rf {\rfloor}

\begin{document}


\author{O. KHORUNZHIY{\footnote{khorunjy@math.uvsq.fr}}\\
Laboratoire de Math\'ematiques\\
Universit\'e de Versailles -- Saint-Quentin\\ 
Versailles 78035, 
\textsc{ France}}

\title{LIMIT THEOREM FOR SUMS \\OF RANDOM PRODUCTS
 }
\maketitle

\begin{abstract}
We study asymptotic behavior of  the sums $Z_{N,m}(\beta)$ 
determined by the formula
$$
Z_{N,m}(\beta) = 
\sum_{i=1}^N 
\prod_{j=1}^{\lambda_i^{(m)}} \exp\{ \beta \eta_{i,j}\}\, \, ,
$$
where $\{\lambda_i^{(m)}\}$ and $\{\eta_{i,j}\}$ are jointly independent
random variables such that their laws do not depend on the subscripts.
We assume that
$\lambda_i^{(m)}$ take values
$ \{0,1,2,\dots\}$ and the mathematical expectations are 
${\E} \lambda^{(m)} =m$ and
$\E \eta = 0$.
The random variable $f^{(\beta)}(N,m) = - (\beta m)^{-1} \log
Z_{N,m}(\beta)$ can be regarded as the free energy of one more 
generalization of the Random Energy Model (REM) widely known in  statistical
mechanics. Imposing natural restrictions on the laws of   $\l$
and $\eta$, we prove existence of the limit  
{\mbox {$f^{(\beta)} = \lim_{m\to\infty} f^{(\beta)}(2^{m},m)$}}
and obtain explicit expression for $f^{(\beta)}$. We study  properties of
$f^{(\beta)}$
 and compare our results
with those obtained for REM and the Erd\H os-R\'enyi limit theorem.

\end{abstract}

\section{Introduction}

The Random Energy Model (REM) was proposed by B. Derrida \cite{D} 
as a simplified version of the spin glasses 
(here the main reference is  \cite{MPV}). 
The statistical
sum of REM can be represented by expression 
$$
Z_{N,n}(\beta) = \sum_{i=1}^N \exp\{{\beta \sum_{j=1}^n \eta_{i,j}}\} = 
\sum_{i=1}^N 
\prod_{j=1}^{n} \exp\{ \beta \eta_{i,j}\}\, \, ,
\eqno (1.1)
$$
where $\{\eta_{i,j}\}$ are i.i.d. random variables with zero mean value 
$\E \eta = 0$. Under assumption that $\E e^{t \eta}$ exists,
it was shown that in the limit $n\to\infty$, $N= 2^n$ 
there is a critical value  $\b_c$ that separates 
two different regimes of the
asymptotic behavior of the free energy 
$f^{(\b)}(N,n)= -{1\over \b n} \log Z_{N,n}(\b)$ \cite{D}.

Almost at the same time, another model was
considered \cite{LK}, where (1.1)  was  assigned a meaning
completely different from that of \cite{D}.
In \cite{LK} $Z_{N,n}$
represents the transparency of disordered layer
in quasi one-dimensional case of quantum particles.
Such a layer is composed of $N$ threads of the length $n$, each thread contains
a random number of obstacles.
The particle
goes along the thread and in the case of low enough energy 
of the particle, the probability
of passing through the thread decreases exponentially with the distance.
However, if one takes $N$ exponentially large with respect to $n$,
one can see some of the particles hitting the other side of the layer.
The number of such particles depends on  the probability distribution
of the obstacles.

In papers \cite{LGP,LK} this model has been rather deeply explored.
Several asymptotic expressions were obtained for various limiting
approximations.
In particular,
the decrement of transparency $f^{(1)}(N,n)$ 
was studied in \cite{LGP}
with $N=e^{qn}$.
Assuming that the probability of large deviations of the sum
$\sum_{j} \eta_{i,j}$  is
 determined by a function
$\phi$, the authors of \cite{LGP} derived an explicit expression for the limiting effective
decrement $f^{(1)}(q)$ in terms of $\phi$;
$$
f^{(1)}(q) = \inf_{  \phi(x) \le q} \{ \phi(x) - x\} - q\ .
\eqno (1.2)
$$
More rigorous derivation  of this formula was  given in   \cite{P}.
The resulting conclusion coming from (1.2) 
is that this model  possesses a  critical value  $q_{cr}$, 
where $f^{(1)}(q)$  changes 
its analytical structure. It is easy to see that $q_{cr} = \phi(x_0)$,
where $x_0$ is determined by equality
$\phi'(x_0) = 1$. Certainly,
$q_{cr}$
 is in one-to-one correspondence with 
$\b_c$ of \cite{D}. 

Regarding the model of $N$ threads of the length $n$ \cite{LK}, it is natural to ask
what will happen in the case, when the length of the threads
is distributed at random.
In this paper we study this generalization of the REM.

In Section 2 we describe the model and formulate results.
In Section 3 we prove the main propositions and derive explicit expression 
for the free energy that generalizes (1.2).
Relations between REM and the Erd\H os-R\'enyi limit theorem are considered in Section 4.
In Section 5 we discuss critical values of $q$ in our generalization of REM
and describe some 
other generalizations of REM related to the version studied in present paper.

\section{The model and main results}

Let us consider the family ${\cal N}= \{ \eta_{i,j}, \, i,j \in {\bf N}\}$ of independent 
identically distributed random variables determined on the same probability space $\Omega$.
We assume that
$
\E \eta_{i,j} = 0$ and that there exists integral
$$
\E e^{t\eta_{i,j}} =  e^{\varphi(t)} <\infty 
\quad {\hbox{for all }} \quad t\in I\subset [0,\infty), 
\eqno (2.1)
$$
where $\E = \E^{(\N)} $ denotes the mathematical expectation with respect to the measure
generated by the family ${\cal N}$.
We denote by $\phi(x)$ the Fenchel-Legendre (or Cramer) transform of $\varphi(t)$ \cite{DZ};
$$
\phi(x) = \sup_{t\ge 0} \{ xt - \varphi(t)\}.
$$

Let us also introduce a family $\L = \{ \l_i^{(m)}, \, i,m\in {\bf N}\}$ of independent random
variables also independent from $\N$. We assume that these random variables take
integer non-negative values and their distributions
do not depend on $i$;
$$
\P \{ \l^{(m)}_i = l\} = p^{(m)}_l, \quad l= 0,1,2,\dots \ \ .
$$
 We denote the mathematical expectation $\E^{(\L)}$ by the 
same sign $\E$.
We assume  that 
$
\E \l_i^{(m)} = m$
and there exists such a function $\chi(t)$ that 
$$
\E e ^{t\l_i^{(m)}} = e^{m \chi(t) (1+o(1))}, \,\, {\hbox{ as }} \,\, m\to\infty.
\eqno (2.2)
$$
We denote $\psi(s) = \sup_{t\ge 0} \{st - \chi(t)\}$. 
\vskip 0.7cm 
{\bf Theorem 2.1}. {\it Let us consider the normalized sum
$$
\hat Z_{N,m}(\beta) = 
{1\over N} \sum_{i=1}^N 
\prod_{j=1}^{\l_i^{(m)}} \exp\{ \beta \eta_{i,j}\}\, \, 
\eqno (2.3)
$$
and determine random variable
$$
\hat f^{(\b)}(N,m) = {1\over\b  m} \log \hat Z_{N,m}(\b)\, .
$$
If $N_q= \lf e^{qm}\rf, \, q>0 $ is the integer nearest to $e^{qm}$, then with probability 1
$$
\lim_{m\to \infty} \hat f^{(\b)}(N_q, m) = \hat f^{(\b)}(q),
\eqno (2.4)
$$
where the limiting $\hat f^{(\b)}(q)$ is determined by relations
$$
\hat \fb(q)= {1\over  \b }\ \inf_{\nu(x)\le q} \{ \nu(x) -\b x \},
\eqno (2.5)
$$
and
$$
\nu(x) = \inf_{y\ge 0} \left[ \psi(y) + y \phi\left({x\over y}\right) \right].
\eqno (2.6)
$$
}

\vskip 0.5cm

Let us compare this result with those obtained for REM (1.2).
First it should be  noted that in each of these models the parameters $\beta$ and $q$ 
are almost equivalent. Indeed, if one divides both parts of  (2.5) by $\beta$,
then one obtains (1.2) with $q$ replaced by $\b q$.
This means that to study the properties of the limit (2.4),
we can restrict ourselves with  $\beta = 1$ and consider (cf. (1.2))
$$
\hat f(q) = \inf_{\nu(x)\le q} \{ \nu(x) - x\}\,  .
\eqno (2.7)
$$
However, in (2.5)
we prefer to keep these two parameters present in order to study  
the transition from REM to Erd\H os-R\'enyi partial sums (see Section 4).

The second observation that relates REM and Erd\H os-R\'enyi sums concerns
function  
$\nu(x)$ (2.6). It has already appeared  \cite{K} in the studies
of the stochastic version of the ER limit theorem (see relation (4.5) 
at the end of this paper).
This  transformation  can be regarded as
a convolution  of two functions $\psi(x)$  and $\phi(x)$.
We study its properties in Section 4.
Let us just note here that if one determines $\psi(x) $ as
$$
\psi(x) = \cases{0, & if $x\in (0,1)$, \cr
+\infty , & if  $x\ge 1$,\cr}
\eqno (2.8)
$$
then expression in square brackets in (2.6) is determined for {\mbox{$0\le y\le 1$}}
and in this case $\nu(x) = \phi(x)$.

Finally, let us note that function $\nu(x)$ determines the 
probability of large deviations of the random variable 
$$
R_{i,m}= \sum_{j=1}^{\l_i^{(m)}} \eta_{i,j}.
\eqno (2.9)
$$ 
Representation of $\nu(x)$ different from (2.6) can be found,
for example, in [5, p.52] or \cite{M}.
This reveals once more relations of the REM with the theory of large deviations
\cite{DZ}
and makes hints to some other generalizations of REM (see Section 5).
Certainly, appearance of the critical values $q_{cr}$ in the both of the models
(1.1) and (2.3) is due to the restrictions $\phi(x)\le q$ and $\nu(x)\le q$,
respectively.
This is what makes the difference between REM and 
the standard theory of large deviations.



\section{Proof of Theorem 2.1}

We follow the reasoning of \cite{P} modified along the lines indicated in \cite{K}.
Let us start with the estimate from below.
It is clear that for any $x>0$
$$
\Z \ge e^{\b mx} {1\over N}\sum_{i=1}^N \sum_{l\ge 1} 
\I \{ \l^{(m)}_i= l\} \cdot  \I\left\{ \sum_{j=1}^l \eta_{i,j}\ge
mx\right\}\equiv  
e^{\b mx}\mu(N,m;x),
\eqno (3.1)
$$
where $\I\{A\}$ is the indicator function of the event $A\subset \Omega$.
The estimates of $Z$ from above and below are based on the properties of 
$\mu(x)=\mu(N,m;x)$. So, let us study this variable in more details.

\subsection{Properties of $\mu$}

Let us start with the average value of $\mu$. Regarding independence of 
$\eta$ and $\l$, we obtain equality
$$
\E \mu(N,m;x) = \sum_{l\ge 1} 
\P \{ \l^{(m)}_i= l\} \cdot  \P\left\{ \sum_{j=1}^l \eta_{i,j}\ge
mx\right\} = S(x).
$$

Applying to the last probability 
relation proved by 
Bahadur and Ranga Rao \cite{BR}, we can write that 
$$
S(x) =  
 \sum_{l\ge 1} p^{(m)}_l \exp\{-m{l\over m} \phi(x{m\over l})(1+o(1))\}.
\eqno (3.2)
$$

Let  $y'$ and $\d$ be  such that $\nu(x) = \psi(y') + y'\phi(x/y')$ and  $y'-\d\ge x$.
 The function $y\phi(x/y)$ is monotone and 
we can write inequality
$$
S(x) \ge  \sum_{l\ge m(y'-\d) }  p^{(m)}_l 
\exp\{-m{l\over m} \phi(x{m\over l})(1+o(1))\} \ge
$$
$$
\exp\{-m(y'-\d) \phi({x\over y'-\d})(1+o(1))\}
\sum_{l\ge m(y'-\d)} p^{(m)}_l  
=
$$
$$
\exp\{-m[\psi(y'-\d) + (y'-\d) \phi({x\over y'-\d})(1+o(1))]\}.
$$
Here we have used relation 
$$
\sum_{l\ge mt} p_l^{(m)} = e^{-m \psi(t) (1+o(1))}
$$
that follows from (2.2) by the G\"artner-Ellis theorem (see e.g. \cite{DZ}).
Now it is easy to see that given $\vep>0$, one can find such $\d$
that 
$$
S(x)\ge e^{-m \nu(x) (1+\vep)}.
\eqno (3.3)
$$
From another hand, inequality
$$
S(x)\le e^{-m \nu(x)(1-\vep)}
\eqno (3.4)
$$
also holds. 
This follows from the observation that the sum in (3.2) is of the order
of the sum
$$
\sum_{ m(y'-\d)\le l \le m(y'+\d) }  p^{(m)}_l 
\exp\{-m{l\over m} \phi(x{m\over l})(1+o(1))\} \le
$$
$$
\exp\{-m (y'+\d) \phi({x\over y'+\d})(1+o(1))\}
\sum_{ m(y'-\d)\le l \le m(y'+\d) }  q^{(m)}_l  = 
$$
$$
\exp\{-m [\psi(y'-\d) + (y'+\d) \phi({x\over y'+\d})](1+o(1))\}.
$$
Then (3.4) follows from continuity of the functions $\psi$ and $\phi$.
Therefore, we can write that
$$
S(x) = e^{-m\nu(x)(1+o(1))} \quad {\hbox{as}}\quad m\to\infty.
\eqno (3.5)
$$

We complete this subsection with computation of 
the variance   of $\mu= \mu(N,m;x)$.
It is easy to find that 
$$
\E \mu^2 - [\E \mu]^2 = -{1\over N} S^2 + {1\over N} S=
S^2\left( -{1\over N} + {1\over N S}\right),
\eqno (3.6)
$$
where we denoted $S= S(x)$.

\subsection{Estimates of $\hat Z$}

Turning back to (3.1) and using (3.3), we can write that
$$
{1\over m} \log \hat Z_{N,m}(\b) \ge -\left\{\nu(x) - \b x\right\} + 
{1\over m} \log {\mu(x)\over S(x)}.
\eqno (3.7)
$$
We have parameter $x$ in our disposition to find the most precise estimate.
The computations that follow will show that such an estimate is obtained
when we take the minimal value of the sum $F(x) = \nu(x) - \b x$ over $x$
satisfying condition $\nu(x) \le  q$.
Let us denote by $x'$ the value of $x$ that provides this  minimum.

Since we assume all functions entering (2.5) to be continuous,
then for any given $\vep>0$ there exists such $\t x < x'$ that $x'  - \t x > \d$
for some $\d>0$ and $F(\t x) - F(x') < \vep$. 
Then  we derive from (3.7) that
$$
 {1\over m} \log \hat Z_{N,m}(\b)  \ge  - F(\t x)  + {1\over m} \log {\mu(\t x)\over S(\t x)}.
\eqno (3.8)
$$

Now we turn to the estimate from above.
To do this, let us first show  that if $x'$ is such that $ \nu(x')>q$, then 
with probability 1 all random variables $ \sum_{j=1}^{\l_i} \eta_{i,j}$
are bounded by $mx'$.  Indeed, denoting
$$
A = \bigcap_{i=1}^N \left\{ \omega: \sum_{j=1}^{\l_i} \eta_{i,j}< mx'\right\},
$$
we can write that
$$
\Pr(A) = [ 1 - \Pr\{ \sum_{j=1}^{\l_i} \eta_{i,j}\ge mx\}]^N = [ 1 - S(x')]^N
= [1-e^{-m\nu(x')}]^{mq}.
$$
Then our proposition follows and $\hat Z_{N,m}(\b) \le e^{m\b x'}$ with probability 1.

Let us cover the interval $(0,x'-\d)$, where $\d$ is as in (3.8),   by $M$ intervals
$(a_l,a_{l+1})$ of the length
$r/M$, $ r= \nu(x'-\d)$ with  $a_0= 0, a_M=r$, we can write that 
$$
\hat Z_{N,m}(\b) \le \sum_{l=0}^{M-1} e^{m\b a_l} {1\over N}\sum_{i=1}^N
{\bf I} \left\{ {1\over m}\sum_{l =1 }^{\l_i^{(m)}} \eta_{i,j} \in (a_l,a_{l+1})\right\}\le
$$
$$
 \sum_{l=0}^{M-1} e^{m\b a_l} \mu (N,m; a_l) = 
 \sum_{l=0}^{M-1} e^{m\b a_l}  e^{-m\nu( a_l) (1+o(1))}{\mu( a_l)\over S( a_l)}\le
$$
$$
e^{-m\inf_l [\nu( a_l)-\b  a_l](1+o(1))} \sum_{l=0}^{M-1} {\mu( a_l)\over S( a_l)}\ .
$$
Combinig this inequality with (3.8), we can write that
$$
\vert {1\over m} \log \hat Z_{N,m}(\b)  +  F(\t x) \vert \le  {1\over m} 
\log \left\{M \max_l {\mu(a_l)\over S(a_l)} \right\}.
\eqno (3.9)
$$
The final remark concerns the difference ${\mu(\t x)\over S(\t x) } - 1$
that is small with {\mbox{probability 1.}} This follows from the estimate (3.6) and the observation
that
$N= e^{qm}$ and  $NS (\t x) = m^{q-\nu(\t x)} = O(e^{\d m})$ as $m\to\infty$.
Theorem 2.1 is proved.

\section{Relations with Erd\H os-R\'enyi limit theorem}

Let us rewrite definition (2.3) in the form 
$$
\hat Z_{N,m}(\b) = {1\over N} 
\sum_{i=1}^N e^{\b R_{i,m}}\ ,  
\eqno (4.1)
$$
where $ R_{i,m}$ are given by (2.9). 
The presence of the exponential function in the sum of (4.1) 
gives more weight to  those terms that have maximal value
of $R_{i,m}$. Then one can say that $\hat Z$ (4.1)
represents an  {{\textsl{arithmetic  mean of maximums}} }
of random variables $e^{\b R_i}$.

From another hand, the {{\textsl{maximum of arithmetic mean values}}} of random
variables is known in probability theory \cite{ER}. Namely, 
the sums
$$
R_{i}^{(n)} = \sum_{j=1}^{n} \eta_{i+j}, \quad i= 1,2,\dots, N-n
\eqno (4.2)
$$
are referred to as the partial sums of $N$ random variables $\eta_i$,
and the limit theorem of Erd\H os-R\'enyi
establishes existence of the limit
$$
\lim_{n\to\infty}\ \ \max_{i=1,\dots, N-n} \ 
\left\{ {1\over n} R_{i}^{(n)}\right\} = \gamma_q, \quad 
N= \lf e^{qn}\rf .
\eqno (4.3)
$$ 
The limit $\bar \g = \g_q$ is 
such that $\phi(\bar \g) = q$ \cite{ER}, where 
$\phi$ is determined by $\varphi$ (2.1) and 
one  takes $\eta_{i,j} = \eta_{i+j}$ \cite{ER}.
If one replaces $n$ by $\l_i^{(m)}$ and by $m$ in (4.2) and (4.3), respectively,
then one obtains random variables $R_{i,m}$  that generalize partial sums
(4.2).
A result analogous to  the Erd\H os-R\'enyi limit theorem can be proved for $R_{i,m}$.
Namely, it is shown in \cite{K} that the random variables
$$
X_{N,m} =  \max_{i=1,\dots, N} \ 
\left\{ {1\over m} R_{i,m}\right\}, \quad 
\eqno (4.4)
$$
converges with probability 1 as $m\to\infty$, $N = \lf e^{qm}\rf$  to a limit
$\tilde \g=\g_q(\phi,\psi)$  determined by relation
$$
\inf_{y \ge 0} \left[ \psi(y)+ y \phi\left({\tilde \g \over y}\right)\right] = q.
\eqno (4.5)
$$

Taking into account these results,
it is natural to ask what is the model that interpolates two
random variables: the statistical sum $\hat Z_{N,m}(\b)$  of REM  (4.1) and 
the partial sums  $X_{N,m}$ of ER limit theorem.
It is easy to see  that the random variable
$$
Y_{N,m}^{(k)}(\beta) = 
\left[{1\over N} \sum_{i=1}^N 
 \left( e^{ \beta R_{i,m}}\right) ^{k\over m} \right]^{1/k}\ 
\eqno (4.6)
$$
represents of the possible interpolations needed.

Indeed, if $k=m$, then  the random variable
$
\log Y^{(m)}_{N,m}(\b)  = -\b f^{(\b)}(N,m)
$
converges to $ -\b f^{(\b)}(q)$ (2.4) as $m\to\infty$.
In the case of  $k>>m$ one can expect $Y_{N,m}^{(\b)}$ to be 
close to $\exp\{X_{N,m}(q)\}$. In the general case,
one can prove the following proposition.

\vskip 0.5cm
{\bf Theorem 4.1} 

{\it If one considers (4.6) in the limit $k,m\to\infty$ such  that $k/m\to\a, k\ge m$, then 
under conditions of Theorem 2.1
}
$$
\lim_{m} \log Y_{N,m}(k)(\b) = 
-\inf_{\nu(x)\le q} \left\{ {\nu(x)\over \a} - \b x\right\} \ ,
\eqno (4.7)
$$
{\it where $\nu(x)$ is determined by relation (2.6).
If $k\gg m$, then}
$$
\lim_{m\to\infty} \log Y_{N,m}^{(k)}(\b) = -\inf_{\nu(x)\le q}\{-\b x\} 
= 
\b\nu^{-1}(q).
\eqno (4.8)
$$

\vskip 0.5cm
\noindent We do not present the proof of  Theorem 4.1
because it does not differ much from the proof of Theorem 2.1.
 
\vskip 0.5cm 
The meaning of this theorem is very simple. When introducing $k$ into (4.6),
one changes the value of $\b$ in expression 
$1/\b \log \hat Z_{N,m}(\b )$. From this point of view, the second limiting transition
of Theorem 4.1 $k/m\to\infty$ corresponds to the limiting transition of (2.4)
with $m\to\infty, \beta\to \infty$. Certainly, this limit
leads one to the maximum of partial sums of the Erd\H os-R\'enyi type; that  is
 proved in {\mbox{Theorem 4.1.}}

The Erd\H os-R\'enyi limit theorem is studied in many aspects (see, for example,
\cite{DD,DDL,N} and references therein) and in a series of applications
\cite{C,RV}. Its relation  with the REM can
be useful in further studies.

Let us study properties of  function  $\nu(x)$ (2.6).
To do this, it is convenient to consider the simplest case
$$
\psi(x) = \cases{ 0, & if $x\in (0,1)$ ,\cr
x(\log x -1) +1, & if $ x\in [1,\infty)$\cr}
$$
that corresponds to the Poissonian distribution of $\l^{(m)}$.
It is easy to  see that $\nu(x)$ is positive,   monotone increasing function such that 
$\lim_{x\to 0} \nu(x) = 0$. 

Also, it is not hard to show that $\nu(x) \le \phi(x)$ and $\tilde \g = \nu^{-1}(q) \ge
\phi^{-1}(q) = \bar \g $. This inequality can be easily explained by the observation that 
the presence of additional "randomness" in the sums (4.2) with respect to the sums
(4.3) increases the value of the maximum in the Erd\H os-R\'enyi 
limit theorem.  In this connection,
it is interesting to study the difference between the critical values of $q$
determined by (1.2) and (2.7), respectively.

\section{Discussion}
We have studied a version of the Random Energy Model
\cite{D,LGP,LK}.
The generalization of (1.1)  is to replace the given number $n$ of factors
by a random variable $\l_i$. Let us describe here how this model
arises from the original arguments presented in \cite{CD}.
There (1.1) has been considered as the statistical sum  of the polymer chains 
determined as  branches of the 
 Cayley
tree.
The energy of a chain is given by the sum of random weights along the
corresponding branch.
 If the tree has the cardinality number $K$, then
it is natural to consider (1.1) with $N= K^n$.
A version of 
the sparse random matrices has been  proposed to describe
this model. In this model, one considers the product of $n$
matrices $B_i$ of dimensions $N\times N$;
the matrices are such that there are exactly $K$ non-zero elements
in each line, and the positions of non-zero elements are chosen at random
at each $B_i$.

It should  be noted  that
in random graph theory, 
another ensemble of sparse random matrices 
arises. Namely, the most studied model of 
the random graph is such that in a graph with $N$ vertices \cite{BB},
 each edge is present with probability $m/N$
or absent with probability $1-m/N$.
In this case  
the adjacency matrix $A$ of the graph is 
a real symmetric matrix 
of dimension $N\times N$
and each its element above the diagonal
is non-zero with probability $m/N$ and zero with probability
$1-m/N$; 
then $(A)_{ij} = a_{ij}, i\le j$ are independent random variables such that 
$$
a_{ij}^{(m,N)} = \cases{1, & with probability ${m\over N}$\ ,\cr
0, & with probability $1-{m\over N}$\ .\cr}
\eqno (5.1)
$$
Obviously, in such a matrix  the average number of non-zero elements per line
is $m$.

When regarding the model of polymer chains with matrices $B$
replaced by $A$, 
one can easily see that 
the length of the polymer chain is given by a random variable.
Thus, we obtain a generalization of REM (1.1),
where the parameter $n$  is replaced by  random variables $\l_i$ 
with mathematical expectation
$m$.  Since the fixed value of $n$ is replaced by a random variable,
the fluctuations of the energy are
greater than that of the REM.
This implies inequality $\nu(x) \le \phi(x)$  mentioned at the end
of Section 4. 
Then the limiting value of the maximum in the Erd\H os-R\'enyi limit theorem
increases.

It is natural to expect the critical value $\hat q_{cr}$ 
determined by (2.7) to be greater than $q_{cr}$.
This conjecture can be supported by some particular cases of $\phi$ and $\psi$ 
with large values of $\beta$. However, $\hat q_{cr}$ in the general case,
as well as the properties of the convolution $\nu(x)$
are to be studied.

Let us discuss two more generalizations of (1.1).
The first one resembles the sparse random matrices \cite{K1}
that also can be considered as a weighted adjacency matrices
of a random graphs. One of the possible model here is obtained when the random variables
$R_{i,m}$ (2.9) are replaced by
$$
R^{(1)}_{i,m} = \sum_{j=1}^m a_{i,j}^{(m,N)} \, \eta_{i,j} \ .
$$
This generalization can be treated by the same
way as it is done in Section 3. More interesting
situation arises when one considers
$$
R^{(2)}_{i,m} = \sum_{j=1}^m a_{i,j}^{(m,N)} \, \eta_{j} \ ,
$$
where $\{\eta_j\}_{j\in {\bf N}}$ is a family of independent random variables (4.2).
However, in this case $R_i$ are correlated that makes
analysis of corresponding free energy
more difficult.

As it is  pointed out, in the case of jointly independent random variables $R_{i,m}$,
the limiting expression (2.7) is determined by the function $\nu(x)$ controlling
probability of large deviations of $R_{i,m}$. 
Therefore, if one replaces $R_{i,m}$
by a family of $N$ independent random processes
${\cal R}_{i}(t)$, the crucial role in the limit
$t\to\infty$, $N= \lf e^{qt} \rf$  will be played by their large values.
It would be interesting to derive expressions analogous to (2.7) and determine
the critical values of $q$ in this model.

At last, let us repeat that we have traced out relations between the limit 
theorem of REM with that of Erd\H os-R\'enyi. In fact, the corresponding result 
(see Theorem 4.1) 
is familiar 
and almost evident to physicists. This is because the case when we obtain Erd\H os-R\'enyi
expression from (4.1)
corresponds to the limiting transition $\b\to\infty$ that immediately
leads one to the maximum of partial sums of random variables.
Up to our knowledge,  this correspondence has not been stated before
in explicit form. The sums
$\sum_{i=1}^{N-n} e^{\b R_i^{(n)}}$ have been  considered in \cite{T} in connection with 
Erd\H os-R\'enyi limit theorem,  but the limit $\beta\to\infty$ has not been regarded there.


\end{document}